\documentclass{preprint}

\usepackage{tikz-cd}				
\usepackage{frege}				    
\usepackage{textgreek}				
\usepackage{yfonts}					
\usepackage{endnotes}
\usepackage{fge}				    
\usepackage{enumerate}
\usepackage{cleveref}

\usepackage{parskip}

\newcommand{\hypjdg}[2]{#2 \; (#1)}

\title{Frege's theory of types}

\journal{Manuscrito}

\volume{46}

\issue{(4)}

\pages{1--47}

\when{2023}

\doi{10.1590/0100-6045.2023.v46n4.bb}

\author{Bruno Bentzen}

\department{School of Philosophy}

\university{Zhejiang University}

\country{Hangzhou, China}

\email{bbentzen@zju.edu.cn}

\newcommand{\funtype}[2]{#1\to#2}
\newcommand{\infuntype}[2]{(#1\to#2)}

\begin{document}

\nocite{frege1962grundgesetze}

\maketitle

\begin{abstract}
	It is often claimed that the theory of function levels proposed by
	Frege in \textit{Grundgesetze der Arithmetik} anticipates the hierarchy of types that
	underlies Church's simple theory of types. This claim roughly states that 
	Frege presupposes a type of functions
	in the sense of simple type theory in the expository language of~\gga. However, this view makes it hard 
	to accommodate function names of two arguments and view functions as incomplete entities. I propose and defend an alternative interpretation of 
	first-level function names in \gga~into simple type-theoretic open terms 
	rather than into closed terms of a function~type. This interpretation offers a still unhistorical but more faithful type-theoretic approximation of Frege's theory of levels and can be naturally extended to accommodate second-level functions. It 
	is made possible by two key observations that Frege's Roman markers 
	behave essentially like open terms and that Frege lacks a
	clear criterion for distinguishing between Roman markers and function names. 		
\end{abstract}

\section{Introduction}\label{intro}

It is often claimed that the theory of function levels introduced by \cite{frege1962grundgesetze} in his two-volume magnum opus \textit{Grundgesetze der Arithmetik} anticipates the hierarchy of types that underlies the simple theory of types developed by \cite{church1940formulation}.\footnote{
	The main proponents of this view are \citet[p.148]{quine1955frege}, \citet[p.330]{resnik1965frege}, \citet[pp.44--50]{dummett1973frege}, \citet[pp.8--9]{martinlof1993leiden}, and~\citet[ch.2,\S1.1]{klev2014categories}. It should be noted, however, that an earlier anticipation is found in Schr\"oder, as \citet[407]{church1939schroeder} and \citet[p.148]{quine1955frege} observe. \citet[\S5.4]{landini1998russell} proposes a formulation of the formal system of \textit{Principia Mathematica} as a simple theory of types.

} 
The claim is roughly that a ground type of objects and a type of functions is present in the ideography. I believe the anticipation is articulated more clearly in the words of Martin-L\"of:

\begin{quote}
	The first type structure in the modern sense was introduced by Frege. [...] 
	He
	did not talk about types, of course, but using modern notation and terminology,
	we would say that he had a type structure, namely the type structure which is
	generated by the following clause:
	If $\alpha_1 , ... , \alpha_n$ are types, then $(\alpha_1 , . . . , \alpha_n)$ is a type.
	In this clause, $n$ is allowed to be $0, 1$, etc. In particular, for $n = 0$, we get the type
	$( )$, which is Frege's type of objects. Then putting $( )$ for each of $\alpha_1 , ... , \alpha_n$, we
	get a type of the form $(( ), ... , ( ))$, which is the type of n-ary functions taking n
	objects into an object again. This can be repeated, of course, arbitrarily high up,
	so we get, for instance, the important type $((( )))$. A function of this type takes a
	unary function from objects to objects into an object, so it is the type of Frege's
	course-of-values operator, the Wertverlauf operator.~\cite[p.8]{martinlof1993leiden}
\end{quote}

It is instructive to contrast those two theories as follows. 
In Frege's theory of function levels, each function of one or two arguments is schematically categorized into first, second, and third levels depending on whether its domain of argument is strictly restricted to that of objects, first-level functions, or second-level functions, respectively. 
On the other hand, in Church's simple type theory only unary functions are formally part of the hierarchy but the formation of functions of multiple arguments is made possible by allowing functions to be values of other functions, in contrast with Frege. More importantly, higher functions are formed through iterations of the function type $\funtype{\sigma}{\tau}$ when $\sigma$ itself is a function type. 
The claim above therefore implies that any function name in the expository language of \gga~can be interpreted as a term of a function type in the sense of simple type theory. For the lack of a better name I shall call this view the `function type interpretation'. 

The main opponent of the function type interpretation is \cite{church1939schroeder}. 
According to him, the view that Frege has a type of functions is based on a misunderstanding. 
Church objects that Frege's view of functions as incompleted abstractions that require completion prevents them from being regarded as actual objects or completed abstractions: 

\begin{quote}
	With Frege a function is not properly an (abstract) object at all, but is a sort of incompleted abstraction. A function is \textit{unges\"attigt} \textendash~is like an incomplete symbol in that it requires something additional, argument or quantifier, to complete its meaning \textendash~but nevertheless partakes sufficiently of the nature of an object to be represented by a variable. The
	notion of a function as an actual object, a completed abstraction in
	intension, is not contemplated by Frege; [...] \cite[p.408]{church1939schroeder}
\end{quote}

We may thus say with \citet[p.74]{klev2014categories} that Church opposes the view that a function in the ideography is an object of a type in the sense of simple type theory. However, Klev dismisses Church's accusations by noting that Frege allows his functions be the arguments to other functions and in the scope of quantification, concluding that it would be hard to deny that the functions that constitute Frege's theory of function levels can form a type. 

Building on Church's criticism of the function type interpretation, 
my aim in this paper is to propose and defend an alternative interpretation of Frege's theory of function levels into simple type theory in which first-level function names are viewed as open terms in a context depending on free variables rather than as closed terms of a~function~type. 
This interpretation can be naturally extended to second-level functions following Church's formal account of abstraction in terms of rules of inference. Thus, I respond to Klev's challenge concerning the place of higher functions and quantification by interpreting second-level function names as inference rules from open terms as premises to closed terms as conclusions.

In other words, my proposed interpretation, which I call the `open term interpretation', is developed without a type of functions, as opposed to the widespread function type interpretation of Frege's theory of function levels. 
Before I put forward the open term interpretation as a more accurate way of understanding Frege's theory of function levels in the simple type-theoretic setting of~\cite{church1940formulation}, I will argue that the function type interpretation is problematic for two main reasons: it cannot properly accommodate function names of two arguments and it goes against Frege's firm conviction that functions are incomplete entities that are in need of completion. This is the negative part of the paper. My positive contribution is that to substantiate the open term interpretation we must delve deeper into Frege's accounts of judgment, functionality, and generality in \gga. The following are two important theses that we will establish along the way: 

\begin{enumerate}[(i)]
	\item Roman markers in the ideography behave like open terms in simple type theory;
	\item Frege lacks a criterion for distinguishing between Roman markers and function~names. 
\end{enumerate}
 
It is important to note that I am not proposing this type-theoretic rendering of the ideography as a portrayal of Frege's actual views on functions. This is certainly not the case. Still, as it draws out the utmost implications of his conviction that function names are incomplete expressions and sheds light on his confusion between Roman markers and function names, I believe this interpretation should be of interest in its own right. 
With the open term interpretation we will solve our initial worries that functions of two arguments and Frege's view of functions as incomplete entities cannot be type-theoretically represented. 
While it has some difficulties dealing with third-level functions, I believe it remains the most faithful approximation of a type-theoretic interpretation of Frege's theory of levels. 

Against my undertaking in this paper it may be correctly observed that no function name actually occurs in the object language of the ideography. Instead, function names appear only in the expository language of \gga, while the language of the formal system only contains names of values of functions, which are names for objects.\footnote{
	See e.g. \citet[p.19]{landini2012frege} and \citet[pp.16-17]{cadetpanza2015logical}. Landini sees the expository language of \gga~as a form of metalanguage since for him small Greek letters are metalinguistic variables. It is unclear to me whether Frege had a full grasp of the difference between object language and metalanguage in the same way that we do today. Perhaps the notion of metalanguage is even intelligible to him, given that, as Cadet and Panza note, it goes against his universalist conception of logic. Either way, the fact remains that the  expository language of \gga~should be distinguished from its object language. For Landini, adopting metalanguage instead of object language does not contradict this universalist concept of logic. } 
The objection is that when `$\Phi(\Delta)$' appears as a name of an object in the ideography, for example, we cannot say that the expression `$\Phi$' appears as a name of a function in Frege's object language. For him only an expression requiring completion such as `$\Phi(\xi)$' can be a function name. Frege is explicit about this issue in a letter to Russell from 13 November 1904 to be discussed in~\Cref{objections}. 
In sum, the problem in question here is that these incomplete names can never occur in the object language of the ideography since they require completion. 
This observation does not undermine the investigation of Frege's theory of function levels within the context of simple type theory that I aim to undertake in this paper. Frege's theory of function levels clearly falls outside the object language of \gga, since it is about when his function names can be properly filled, and therefore so does my investigation of its nature. 
It should be noted that when I speak of the function names of the ideography in what follows, I am referring specifically to the names of functions assumed in the theory of function levels that operates behind the scenes in \gga. 
For example, if `$\Phi(\Delta)$' and `$\Delta$' are names in the object language, then `$\Phi(\xi)$' must be a function name in the expository language. 
The aim of this paper is to investigate how, if possible, the structure of Frege's theory of function levels can be analyzed within the framework of simple type theory. 
Still, one of the outcomes of my investigation is that Frege's design choice of admitting Roman markers in the object language~(as we shall soon see in~\Cref{romanmarkers}) but keeping function names only in the expository language is unjustified given that his views in \gga~are deficient in providing a criterion for distinguishing Roman markers and function names.

The rest of this paper is structured as follows. 
In \Cref{typetheory}, I introduce some basic concepts from simple type theory, but only loosely following Church's original formulation. Instead my exposition borrows from the presentation of type theory in the tradition of \cite{martinlof1975intuitionistic} and focuses on the role of open terms and his distinction between functions in the old-fashioned and modern sense. 
I discuss in \Cref{objectandfunction} Frege's views on functionality, his theory of function levels, and the function type interpretation, closing this section with four objections to the function type interpretation. 
To set the stage for the open term interpretation, it will be necessary to establish the theses (i) and (ii) above first. 
\Cref{romanmarkers} is primarily a defense of (i) based on Frege's account of judgment and his use of Roman markers. \Cref{openterms} supports (ii) with an examination of Frege's views on functionality and generality. 
In~\Cref{simpletypetheory} I put together all the elements of our defense of the open term interpretation and explain it in full for first-level and second-level functions, while also addressing its interpretative problem with respect to third-level functions. 
\Cref{conclusion} offers some concluding remarks. 

\section{Open terms and functions in simple type theory}\label{typetheory}

We begin with a brief overview of simple type theory. My exposition shall be based on the tradition of \cite{martinlof1975intuitionistic} with minor differences in terminology and notation. 
Simple type theory is a formal system restricting the operations of the untyped lambda calculus with the adoption of a type system with some ground types and a type former~$\funtype{\sigma}{\tau}$ for functions that take arguments of type~$\sigma$ and result in values of type~$\tau$. 
We write the typing judgment that states that $a$ is a term of type $\sigma$ as `$a : \sigma$'. 
Saying that $a : \sigma$ is a judgment simply means that it may carry assertive force and be subject to rules of inference. When such a judgment is asserted it is written preceded by a turnstile symbol following Frege's notation: 

\vspace*{-1.3em}

$$\vdash a : \sigma. $$

\vspace*{-1.em}

As in modern logic, the turnstile symbol is not part of the object language of the theory. 
I will say that a judgment $a : \sigma$ is categorical when it is stated without any hypotheses. For the most part, however, a typing judgment takes a general hypothetical form, where a typed term depends on one or more typing assumptions on free variables. In this case it will be written as $\hypjdg{\Gamma}{a : \sigma}$, where $\Gamma$ is a finite list of typed variables $x_i : \sigma_i$ called the context of the general hypothetical judgment. Notice that when the context of a general hypothetical judgment is empty it just boils down to a categorical judgment.

We also say that a typed term $a : \sigma$ is open when it occurs in a general hypothetical judgment and may depend on some of the variables declared in the context. A typed term occurring categorically is said to be closed, meaning only bound occurrences of variables. 
A term is open if it is not closed, so some terms may be called open regardless of whether the variables declared in the context actually occur in them. In fact, contexts can be extended freely and a term $a : \sigma$ formed under a specific context may also be said to depend on an additional typed variable. This property is known as the structural rule of weakening: \\

\vspace*{-1.0em}

\begin{equation} \label{eq:weak} 
\frac{\vdash \hypjdg{\Gamma, \Delta}{a : \sigma}}{\vdash \hypjdg{\Gamma, x : \tau, \Delta}{a : \sigma}.}
\end{equation} 

\vspace*{-1.0em}

In the presence of an open term $b : \tau$ depending on a free variable $x : \sigma$, we may obtain a new closed term, when given a closed $a : \sigma$, by substituting $a$ for $x$ in $b$. This idea is generalized for arbitrary contexts with the substitution principle: \\

\vspace*{-1.em}

\begin{equation} \label{eq:ttsubst} 
\frac{\vdash \hypjdg{\Gamma,x : \sigma, \Delta}{b : \tau} \quad \hypjdg{\Gamma}{\vdash a : \sigma}}{\vdash \hypjdg{\Gamma,\Delta}{b[a/x] : \tau}.}
\end{equation} 

\vspace*{-1.em}

Open terms have no meaning in isolation in the meaning explanations proposed by \cite{martinlof1982constructive} as the intended interpretation of type theory. An open term is said to have type $\sigma$ when it results in closed terms of type $\sigma$ after all the typed variables occurring in the context of the general hypothetical judgment are substituted with closed terms of the same type. 

\subsection{The function type}

The definition of the function type $\funtype{\sigma}{\tau}$ involves a well-established collection of rules of inference that we do not need to dwell on in this paper. The reader interested in a comprehensive exposition is referred to \citet[\S19.6]{nordstrom1990programming}. 
Here it will suffice to note that lambda abstraction, the rule that regulates the construction of functions as lambda terms $\fun x a$, acts on an open term $a$ depending on a distinguished free variable $x$: 

\begin{equation} \label{eq:lambdaabst} 
\frac{\vdash \hypjdg{\Gamma, x : \sigma}{a : \tau}}{\vdash \hypjdg{\Gamma}{\fun{x}a : \funtype{\sigma}{\tau}}.} 
\end{equation}

As pointed out by \citet[p.114]{quine1940church}, this technique is essentially a revival of Frege's distinctive approach to abstraction using value-ranges. 
In a dictionary entry on abstraction, \citet{church1942abstraction} claims that the smooth-breathing notation $\spirituslenis{\varepsilon}f(\varepsilon)$ introduced by Frege to describe value-ranges is one of the precursors of the lambda-notation $\fun{x}a$ he uses for the function terms that inhabit the function type. 
This is hardly a coincidence, given that, according to \citet[p.408]{church1939schroeder}, who, again, denies that Frege has a type of functions, it is the value-range that ``corresponds to the notion of a function as used in mathematics''.

Another point worth contrasting with Frege, to be explored in \Cref{objectandfunction}, is that although the function type only gives us functions of one argument, functions of multiple arities can be defined in simple type theory through multiple unary functions by means of ``Currying'', a technique of dispensing with functions of multiple arguments by allowing functions to have other functions as values~\citep[\S2]{schonfinkel1924bausteine}. For example, the type of binary functions with domains $\sigma$ and $\tau$ and codomain $\rho$ is $\funtype{\sigma}{(\funtype{\tau}{\rho})}$, meaning that to apply such a function one has to first apply it to a term of $\sigma$ and only then to a term of $\tau$. 

\subsection{Functions in the modern and old-fashioned sense}

\citet[pp.37--38]{martinlof1993leiden} contrasts two notions of function in type theory, noting that functions in the old-fashioned sense are given by open terms in a type $\hypjdg{x : \sigma}{a : \tau}$ and functions in the modern sense given by closed terms in a function type $f : \sigma \to \tau$.\footnote{
	This distinction has been discussed in \citet[\S4]{klev2019name} and more recently in \cite{klev202Xspiritus} in the context of the argument move and removal rules proposed by \citet[pp.59--64]{martinlof1993leiden} as alternatives to abstraction and application as inference rules for the function type. } 
The former notion of function is called so because it is related to an old practice, introduced in the 18th century by Bernoulli and Euler, and followed by Lagrange, Fourier, and Boole, among others, of incorporating the argument places in the body of the functional expression. 
The reader is directed to \cite{ruthing1984some} for a useful chronological compilation of some definitions of the concept of function from Bernoulli to Bourbaki and to~\citet[pp.50--53]{martinlof1993leiden} for an extensive historical account of the concept of function. 
From Martin-L\"of's perspective, functions in the old-fashioned sense are characterized by their dependence on arguments, examples of which are the functions given as $sin(x)$ or $(2 + 3x^2)x$. 

In contrast to the old-fashioned practice, today we often define $f(x)=(2 + 3x^2)$ and then refer to $f$ as a function without its arguments, for example. Functions in the modern sense are characterized by the absence of arguments. 
\citet[p.52]{martinlof1993leiden} attributes the modern notion of function to Dedekind and his definition of a function $f$ as a law that transforms an object $a$ into an object $f (a)$, in which a function is viewed as a self-subsistent mathematical operation that transforms every object into a value. It should be stressed that this modern concept of functions is implicit in Church's account of functions as ``completed abstractions''. In the modern sense we say that $sin$ or $\fun x (2 + 3x^2)x$ are functions instead. 
Martin-Löf also adds that it is not clear whether Frege can be said to pursue the old-fashioned or modern notion of function in \gga~\citep[p.53]{martinlof1993leiden}. Since, as we will see in the next section, Frege views in \gga~\S1 a function expression as the part that remains when the argument place is not present in the expression, Martin-Löf maintains that Frege was actually aiming at the modern definition of a function. But he acknowledges that Frege writes function names in the old-fashioned manner, thus concluding that Frege did not have a satisfactory way of denoting functions in the modern sense:

\begin{quote}
	When Frege
	wants to have an expression for such a function, he does not do what we now are
	doing, say define $f (x, y) = x^2 + 3y$ and then speak of $f$ as the function in the
	modern sense. This is the natural thing to do now, but Frege does not have a
	satisfactory notation for functions in the modern sense. The way he handles it is
	by introducing these special Greek symbols, $\xi$ and $\zeta$, and he takes $\xi^2 + 3\zeta$ to be the
	expression of the function. Now of course he has after all put these Greek letters
	into the argument places, so it does not have holes in it any longer, it just looks
	like the usual expression for an old-fashioned function. \cite[p.53]{martinlof1993leiden}
\end{quote}

In Martin-Löf's terminology, the function type interpretation states that Frege's functions must be understood in the modern sense. That is, a function name `$\Phi(\xi)$' in the theory of levels should be regarded as a closed term of some function type $\tau \to \sigma$. Under the open term interpretation to be defended here, however, Frege's functions must be taken in the old-fashioned sense, meaning that `$\Phi(\xi)$' is an open term of a type $\sigma$ depending on a type $\tau$. 
Before we can decide which interpretation should be rejected and which should be endorsed we must take a closer look at Frege's views on functionality.

\section{Frege's theory of function levels}\label{objectandfunction}

Now let me return to the discussion of Frege's theory of function~levels. 
There was an influential view among Frege's predecessors that a function is an analytical expression distinguished by its dependence on variables. It can be found, for instance, in the writings of Lagrange and Boole, and, as can be seen in~\cite{ruthing1984some}, traced back to Euler. 
In an attempt to avoid the pitfalls of formalism, Frege opposes the idea that functions are expressions already in \gga~\S1. 
Frege tells us that one will be tempted to rather say that a function is the meaning or reference of such an incomplete expression. But according to him the problem with this view is that in `$(2+3x^2)x$' we find a variable `$x$', which does not refer to an object as `$2$' does, for example, but, as he says, ``only indicates one indeterminately''.\footnote{ \label{fn:indication}
	``\textit{nur unbestimmt andeutet}''~(\gga~\S1). } 

The argument places of a function name are from this point on are written with small Greek letters such as `$\xi$' and `$\zeta$' throughout \gga. Therefore, such small Greek letters only indicate indeterminately. \citet[\S2.6]{landini2012frege} maintains that these Greek letters are metalinguistic devices used parametrically to facilitate the rule of substitution. They are therefore not part of the object language of the ideography. 
As it will be seen in \Cref{romanmarkers}, Frege will also speak of Roman letters as signs that only indicate something indeterminately later in \gga~\S8. For Frege, these letters have an entirely different purpose. 

Because the argument places of a function name do not refer, Frege rejects the above proposal and concludes in \gga~\S1 that a function is in essence the connection established between the objects whose names we put for a variable and the objects that then appear as referents of the now complete expression. Assuming the old-fashioned convention that the name of a function must always carry an argument place, Frege writes that: 

\begin{quote}
	The nature of the function reveals itself, rather, in the bond
	that it establishes between the numbers whose signs we put for `$x$' and the numbers
	that then result as the reference of our expression [...] 
	The nature of the \textit{function} lies therefore in that part of the expression that is present without the `$x$'. (\gga~\S1)
\end{quote}

\subsection{Functions of first, second and third level}

The view of functions expressed above is stated with greater precision as the criterion of referentiality proposed for function names in \gga~\S29. Frege states four separate criteria for first-level function names of one and two arguments, and second- and third-level function names of one argument. Again, all these function names only occur as part of the expository language of \gga. According to Frege, a function name is referential if it always results in a referential name when each of its argument places are filled with referential names of the appropriate level. 
First-level function names can only be filled with object names, complete expressions which do not in themselves carry argument places. 
Note that these names of objects actually occur in the formal language. 
Second- and third-level function names must be filled with first- and second-level function names respectively. 

How do the level restrictions work in the ideography? Call a well-formed expression~a~name. 
Frege states in \gga~\S28 that for an expression to be a name, it must either be a primitive, defined, or a compound name. In the ideography, a compound name can only be formed out of a function name by the filling of its argument places with other names in accordance with the level restrictions and formation rules of \gga~\S30. 

In total there are eight primitive names admitted by Frege in \gga, all of which are function names having their own levels. If we ignore for a moment the distinction between functions of one or two arguments and focus on their levels we can describe them as follows. 
The primitive first-level function names are those for the horizontal, negation, implication, identity, and definite article operator, and their object argument places are marked with the small Greek letters `$\xi$' or `$\zeta$'.  Following Frege, I respectively write these names as:

\vspace*{-1.em}

$$`\F[1] \xi\text{'}, \quad `\Fncontent \xi\text{'},  \quad `\Fcontent\Fconditional{\Fcontent\zeta}{\Fcontent\xi}\text{'},   \quad `\xi = \zeta\text{'}, \quad `\backslash\xi\text{'}.$$


The primitive second-level function names are those for the first-order universal quantifier and value-range operator. They are written respectively as follows, where their first-level argument places are indicated with the small Greek letter `$\phi$':

\vspace*{-1.em}

$$ `\Fquant{a} \phi(\mathfrak{a})\text{'}, \quad `\spirituslenis{\varepsilon}\phi(\varepsilon)\text{'}.  $$ 

\vspace*{-0.6em}

The only primitive third-level function name considered is that of the second-order quantifier. Its second-level argument place is marked with the small Greek letter `$\mu$'. The letter `$\beta$', on the other hand, binds the first-level argument place of the second-level function name of one argument used to fill the second-level argument place:

\vspace*{-1.em}

$$ `\Fquant{f} \mu_\beta \mathfrak{f(\beta)}\text{'}.  $$ 

\vspace*{-0.6em}

\newpage

As the ideography lacks primitive object names, it is only possible to form compound or defined object names. The formation of compound object names, in particular, is done by filling the argument places of second- and third-level function names in the appropriate way. This is precisely how names of truth values (sentences) and value-ranges, the only object names considered by Frege in the formal language, are formed in the ideography.

\subsection{The function type interpretation}

The theory of levels does prevent a function from having itself as argument, but the doctrine of value-ranges as objects allows self-application to enter through the back door by means of objects of the kind $f(\spirituslenis{\varepsilon}f(\varepsilon))$.\footnote{\label{fregefirstletterrussell}
	See Frege's letter to Russell of 22 June 1902 \citep[pp.131--133]{frege1980correspondence}. }
\citet[p.148]{quine1955frege} claims that without this peculiar doctrine we cannot reproduce Russell's paradox in the ideography. He concludes that the theory of function levels anticipates, to some degree, the theory of types. 
\citet[p.50]{dummett1973frege} states that the Frege's theory of function levels is essentially Russell's simple theory of types formulated in terms of Frege's notion of incomplete expressions.
\citet[\S2]{resnik1965frege} maintains that the result of Frege's theory of function levels is a simple hierarchy of types for incomplete entities with differing kinds of incompleteness and, to make the idea precise, he proposes the inductive definition for Frege's types seen earlier in \citet[p.8]{martinlof1993leiden}:

\begin{quote}
	To be more precise, we
	define types inductively as follows: 
	$\alpha$ is a type; if $t_1, ..., t_n$ are types, then so is the
	ordered $n$-tuple $\langle t_1, ..., t_n \rangle$ ($n > 0$). The complete entities make up type $\alpha$. An
	incomplete entity belongs to the type $\langle t_1, ..., t_n \rangle$ if and only if it may be completed by
	(and only by) entities of types $t_1, ..., t_n$ taken in that order. Thus, a function of two
	arguments of type $\langle \alpha, \langle \alpha \rangle \rangle$, for example, can only take as its first argument a member of
	type $\alpha$ (an object) and as its second argument, a member of type $\langle \alpha \rangle$ (a function). \citep[p.330]{resnik1965frege}
\end{quote}

More precisely, the claim is that the theory of function levels amounts to a hierarchy of simple types formulated in the old-fashioned kind, that is, in terms of Frege's view of functions as ``incompleted abstractions'', to borrow Church's terminology. 
When spelled out in detail the idea amounts to a function type interpretation. It states that the names occurring in \gga~at the object- or meta-level can be understood in the setting of a simple type theory with a ground type of individuals $\iota$ and a function type $\funtype{\sigma}{\tau}$, with the following hierarchy of types reflecting the structure inherent to Frege's theory of function levels: 

\begin{itemize} 
	\item Object names have type $\iota$;
	\item First-level function names have type $\funtype{\iota}{\iota}$; 
	\item Second-level function names have type $\funtype{\infuntype{\iota}{\iota}}{\iota}$; 
	\item Third-level function names have type $\funtype{\infuntype{\infuntype{\iota}{\iota}}{\iota}}{\iota}$. 
\end{itemize}

The eight primitive names of the ideography described above, for instance, are interpreted as constant terms of the function type corresponding to their level. 
Thus, for example, we have, as functions of first-, second-, and third-level respectively, the terms

$$\F \xi : \funtype{\iota}{\iota}, \quad ~ \Fn \xi : \funtype{\iota}{\iota}, \quad~ \backslash\xi :  \funtype{\iota}{\iota}, \quad~ \text{and}$$

$$\Fquant{a}\phi(\mathfrak{a}) : \funtype{\infuntype{\iota}{\iota}}{\iota}, \quad ~ \spirituslenis{\varepsilon}\phi(\varepsilon) : \funtype{\infuntype{\iota}{\iota}}{\iota}, \quad~ \text{and}$$

$$\Fquant{f} \mu_\beta \mathfrak{f(\beta)} : \funtype{\infuntype{\infuntype{\iota}{\iota}}{\iota}}{\iota}.$$

\subsection{Objections to the function type interpretation} \label{objections}

Before I proceed to the defense of the open term interpretation in the next section, I wish to raise four objections to the function type interpretation stated above. The first two concern functions of two arguments and the challenges the function type interpretation faces in accommodating them. The other two concern the modern view of functions and how it contradicts Frege's belief that functions are incomplete entities that require completion.

1. I deliberately left the hierarchy of types above incomplete. I chose to only include functions of one argument because, as we shall now see, the function type interpretation fails to capture Frege's functions of two arguments properly. 
One may wonder if we could not simply assign these functions the type $(\sigma \times \tau ) \to \rho$, where~$\sigma \times \tau$ is the product type which is inhabited by ordered pairs of terms of types $\sigma$ and $\tau$. This seems to be implicit in the suggestion made by \cite{resnik1965frege} and \cite{martinlof1993leiden} in their inductive definition of Frege's types. One problem with this approach is that it does not reflect Frege's process of double completion. 

When introducing functions of two arguments in \gga~\S4, Frege explains that they stand in need of double completion given that a function of one argument is obtained after an initial completion and an object is obtained as a value only after a second completion:

\begin{quote}
	So far only functions with a single argument have been talked about; but we
	can easily pass on to \textit{functions with two arguments}. These stand \textit{in need of double
	completion} insofar as a function with one argument is obtained after their completion
	by one argument has been effected. Only after yet another completion do we arrive at
	an object, and this object is then called the \textit{value} of the function for the two arguments. (\gga~\S4)
\end{quote}
 
There is no particular order for the double completion, so Frege has to use different small Greek letters to make each argument place explicit. 
Given a function $f(\xi,\zeta)$ and an appropriate argument $a$, we always have two possible choices: either we apply $a$ to the $\xi$-argument place and obtain $f(a,\zeta)$ or choose the $\zeta$-argument place and instead obtain $f(\xi,a)$. 

However, given a binary function of type $(\sigma \times \tau ) \to \rho$ we can only perform a simultaneous application of both arguments together as a pair.\footnote{
	\cite{simons2019double} notes that Frege does articulate a notion of simultaneous application with ordered pairs in \gga~\S144, but for double value-ranges and not functions of two arguments. } 
It is impossible to provide an initial completion with one argument and then obtain a function of one argument as the result. 
To get around this difficulty, one might then attempt to assign Frege's functions of two arguments the type~$\funtype{\sigma}{\infuntype{\tau}{\rho}}$ of binary Curried functions instead. 
Here one problem is that Curried functions have a fixed order of application while for Frege a function $f(\xi,\zeta)$ can be given an initial completion as either $f(a,\zeta)$ or $f(\xi,a)$.
Type-theoretically, this is an impossible move because the technique of Currying leaves no room for ambiguity in the determination of the order of application. For $f : \funtype{\sigma}{\infuntype{\tau}{\rho}}$, the two applications above may be represented as $f(a)$ and $\fun x (f(x))(a)$, respectively. 
This shows that the type of $f(\xi,\zeta)$ cannot possibly take the Curried form $\funtype{\sigma}{\infuntype{\tau}{\rho}}$, for $\xi$- and $\zeta$-arguments of respective types $\sigma$ and $\tau$. Otherwise we could only doubly complete $f(\xi,\zeta)$ in that order. 

2. A more serious difficulty with the interpretation of Frege's function of two arguments as Curried functions is that the latter have functions as values, as stressed by \cite{potts1979grossest}. Curried functions are allowed to have a value after a first application but, as seen above, the value of a Fregean function of two arguments is only obtained after the second completion. 
So, Frege's functions of two arguments cannot be consistently typed with function types. 

3. 
Frege does not give functions the self-subsistent treatment that is demanded by the notion of function in the modern sense, that is, as an object of the type $\sigma \to \tau$. 
This point is perhaps best made by \citet[p.408]{church1939schroeder} in his criticism of Frege's view of function. 
We can perhaps say in Church's terminology that functions in the modern sense are properly ``(abstract) objects'' or ``completed abstractions'' characterized by their absence of arguments. They are given by closed terms without dependence on free variables. Due to their complete nature they can be objects of a type while functions in the old-fashioned sense cannot. 

This should be taken as an argument against Martin-L\"of conclusion that Frege was actually aiming at a modern definition of function. When Frege states in \gga~\S1 that a function is given by the part of an expression that is present without the argument place he simply means that the function has an incomplete nature. It is not suggested that a function can be regarded as a mathematical entity that is not in need of an argument. 
Thus, it seems to me that Frege adheres more closely to the old-fashioned conception of functions. 
But, as incompleted abstractions, old-fashioned functions cannot be objects of a type. 

4. The function type interpretation implies that first-level function names are interpreted as terms of type $\iota \to \iota$, second-level function names as terms of type $(\iota \to \iota) \to \iota$, and third-level function names terms of type $((\iota \to \iota) \to \iota) \to \iota$. But considering that in type theory terms are either open or closed, there are two possible readings of the interpretation. 

Open terms seem the correct choice, for they arguably embody in type theory the notion of function in the old-fashioned sense embraced by Frege. However, if a function name `$f(\xi)$' is taken as an open term depending on a free variable `$\xi$', its type cannot be $\sigma \to \tau$. 
Take the definite article operator as an example. In this case, $\hypjdg{x : \iota}{\backslash x : \iota \to \iota}$ would amount to an old-fashioned function that has a function of type $\sigma \to \tau$ as value when the argument place $x$ is substituted by a closed term of type~$\iota$. 
This is not what we are looking for. 

Now, if the function type interpretation is read by means of closed terms, we are led to the implausible view that terms with free occurrences of variables like $\backslash x : \iota \to \iota$ are closed. As a result, the function names of \gga~would have to be interpreted in modern style. 
So, if we wish take the function type interpretation seriously, Frege's definite article operator would have to be given by the closed term $\backslash : \iota \to \iota$, for example. 
However, in a letter to Russell dated 13 November 1904, Frege entertains but ultimately rejects the idea of adopting a modern-style notation for functions. The notation proposed is very similar to the function abstraction notation used for value-ranges, except that the variable binder is indicated with a rough breathing mark, while in value-range names smooth breathing marks are used. So, to give an illustration, the proposal is to write a function such as power of two as `$\spiritusasper{\varepsilon}(\varepsilon^2)$' rather than `$\xi^2$' in his standard notation. But Frege says: 

\begin{quote}
	[T]his notation would lead to the same difficulties as my value-range notation and in addition to a new one. For a range of values is supposedly an object and its name a proper name; but `$\spiritusasper{\varepsilon}(\varepsilon^2 = 1)$' would supposedly be a function name which would require completion by a sign following it. `$\spiritusasper{\varepsilon}(\varepsilon^2 = 1)1$' would have the same meaning as `$1^2 = 1$', and accordingly, `$\spiritusasper{\varepsilon}(\varepsilon^2 = 1)\varsupset$' would have to have the same meaning as `$\varsupset^2 = 1$', which, however, would be meaningless. `$\spiritusasper{\varepsilon}(\varepsilon^2 = 1)$' would be defined only in connection with an argument sign following it, and it would nevertheless be used without one; it would be defined as a function sign and used as a proper name, which will not do. \cite[pp.161--162]{frege1980correspondence}.
\end{quote}

It has been argued by \citet[p.20]{klement2003russell}, \citet[p.22]{nomoto2006methodology}, \citet[p.98]{martinlof2006comments} and \citet{klev2019name,klev202Xspiritus} that in the passage above Frege alludes to an abstraction operation that in effect amounts to Church's lambda abstraction (\ref{eq:lambdaabst}). 
Yet, Frege is reluctant to represent the application of a function name by means of juxtaposition because he thinks the lack of argument places in the body of the expression means that the function name can be used as a complete expression. Under Frege's syntactic analysis, it would have to be an object name. 
Frege insists that it is possible to use `$\spiritusasper{\varepsilon}(\varepsilon^2 = 1)$' in isolation as an object name given that it has no occurrences of argument places while `$\xi^2 = 1$' can be viewed as a function name since the expression needs completion. 
So, to sum up, this closed-term based reading of the function type interpretation is unacceptable to Frege. 

Finally, to wrap up the discussion, it should be noted that if by a ground type we mean the type of complete entities, then Frege implies an identification between closed terms and terms of a ground type in his letter. A similar confusion between closed terms and object terms is reflected in Frege's views on the violation of the law of excluded middle, which, as recently discussed by~\cite{vanrie2020did}, throw light on Frege's rejection of Russell's initial proposal of a theory of types. For Frege, there are no closed terms of the function type. 

\section{Roman markers in \gga} \label{romanmarkers}

Since the function type interpretation is untenable, we shall now turn to the development and support of the open term interpretation. 
The basic idea is that a first-level function name such as `$\backslash \xi$' must be interpreted in the old-fashioned sense as $\hypjdg{x : \iota}{\backslash x : \iota}$. This idea requires careful elaboration and I will return to it on \Cref{simpletypetheory}. 
First, we must discuss and overcome what we may see as the main challenge to this interpretation, namely, the fact that the role of open term is already played by the so-called ``Roman markers'' in the ideography. Thus, the open term interpretation conflates both function names and Roman markers into open terms. If function names and Roman markers were to serve different purposes in \gga, this move would be inconsistent with Frege's own views. 
I will be addressing this problem in the remainder of this section and the next. My answer will be that Frege has no criterion to distinguish between function names and Roman markers. 

\subsection{Roman markers as open terms}

To begin with, let me explain why Roman markers behave as open terms in the ideography. 
Roman markers are incomplete expressions depending on Roman letters. They are not names but only expressions that convey generality. 
This straightforward interpretation that the ideography has Roman letters as free variables in the modern syntactic sense is found in \citet[\S 1]{kemp1998propositions}, \citet[p.38]{klement2001frege}, \citet[\S3]{heck2007frege}, and \citet[p.32,\S5.1]{landini2012frege}. 

We will soon see below that Roman markers can only result in names when all their Roman letters are substituted with appropriate names. 
Yet, Frege frequently states his axioms and theorems with Roman markers to emphasize their general applicability. This can be seen in many places, but especially in \gga~\S47 in his summary of the basic laws and in his proof of the derived laws in, for example, in \gga~\S\S49--52.  
This usage is to be expected once we observe that these are actually axiom and theorem schemes and thus cannot concern names of the object language. There is no precise mention of Roman markers when the first three rules of inference of the system are exposed in \gga~\S\S14-16. But it is already clear in \gga~\S17 from Frege's examination of how the Barbara syllogism fits in the ideography that his rules of inference need Roman markers~\citep[pp.38--40]{klement2001frege}. 

The axioms schemes and rules of inferences in simple type theory are, in the same way, formulated by means of open terms and their applicability is not limited to closed terms. 
When we think of the Roman letters of the ideography as typed variables and its names as closed terms, this property can be expressed type-theoretically by the substitution rule (\ref{eq:ttsubst}). 
Likewise, following \cite{martinlof1982constructive}, an open term is said to have a given type when it always results in closed terms of that type after all their typed variables are substituted with closed terms of the same type. 
Formally, the only noticeable difference between Roman markers and open terms is that the former always have occurrences of Roman letters but, due to weakening, free typed variables need not occur in the latter. Despite the converse not being true, every Roman marker in \gga~can be viewed as an open term. 

\subsection{Roman letters and indication} \label{romanletters}

In \gga~\S17 Frege, emphasizes that Roman letters are used to express generality in a judgment and have scope over the entire expression. 
Frege adds in a footnote that the ``use of Roman letters is hereby explained only for the case in which a judgment-stroke occurs.'' We may then say with \cite{landini2012frege} that each occurrence of the same variable is to be regarded as having the same referent once a referent is assigned and this holds throughout judgments of the whole proof. 
Frege stipulates that Roman letters do not refer but merely indicate: 

\begin{quote}
	I shall call \textit{names} only those signs or combinations of signs that refer to something. Roman letters, and combinations of signs in which those occur, are thus not \textit{names} as they merely \textit{indicate}. A combination of signs which contains Roman letters, and which always results in a proper name when every Roman letter is replaced by a name, I will call a \textit{Roman object-marker}. In addition, a combination of signs which contains Roman letters and which always results in a function-name when every Roman letter is replaced by a name, I will call a \textit{Roman function-marker} or \textit{Roman marker} of a function.~(\gga~\S17)
\end{quote}

Here to say that a Roman marker indicates means that it always result in a name that refers to an object or function when all its Roman letters are replaced with appropriate names. 
Now, for example, we can say that the following expression is a Roman object marker because when `$a$' and `$b$' are replaced by object names the resulting names have objects as references: 
$$\text{`}\Fcontent\Fconditional{\Fcontent\Fconditional{\Fcontent a}{\Fcontent b}}{\Fcontent a \qquad.}\text{\!\!\!'}$$ 

This is not the first time Frege speaks of indication as opposed to reference for expressions. 
Recall that in \gga~\S1 Frege stresses that the argument places of a function name do not refer but only ``indicate indeterminately''. The same thing is said of Roman letters in a preliminary discussion on generality in \gga~\S8. In \S17 Frege states that both Roman and Gothic letters do not refer to an object but only indicate it. 
\citet[p.36]{klement2001frege} thinks that to indicate for Frege is to reveal the kind of supplementation that is needed in an incomplete expression, and concludes that Roman and Gothic letters are semantically no different from Frege's small Greek letters. 
I believe this is correct, though for reasons of space I will only focus on Roman and small Greek letters here. 
Frege sees the difference between a generality of an identity in `$\Phi(x)=\Psi(x)$' and an identity in `$2^2 = 4$' as manifested in the presence of a letter `$x$' that ``only indicates indeterminately'' in the former, since  every sign has a determinate reference in the latter. 
Roman letters not only merely indicate, they also indicate indeterminately, just as the small Greek letters do. In \gga~\S1 Frege explains that argument places of a function name indicate indeterminately because ``[f]or different number-signs put in place of `$x$' we generally obtain different references.'' 
This is precisely what it means to say that a Roman marker indicates. 

\newpage

Now, the matter of the presence of Roman markers in the object language of the ideography is a delicate one, as \cite{cadetpanza2015logical} emphasize. Unlike small Greek letters, Roman letters can occur in a judgment in \gga, so while function names are only part of the expository language, Roman markers belong to the object language. Even Roman letters indicating functions are allowed in the object language~\citep[\S\S2.1--6]{landini2012frege}. 
The problem is that Roman markers are not accepted as names in \gga, and, as a result, do not correspond exactly to the terms of a modern formal system. 
In his reconstruction of the ideography as a modern predicate language, Landini avoids this problem by interpreting names as closed terms and allowing Roman markers to be open terms. 
This resembles the type-theoretic reading of Roman markers I am proposing here. 
However, Landini's approach gives rise to a new interpretative challenge as, since his analysis is grounded on a predicate language, he is led to interpret names preceded by Frege's turnstile sign as formulas (i.e. wffs) of a formal system~\citep[\S2.1,\S2.6]{landini2012frege}. 
I will have more to say about the problem with this interpretation after an analysis of Frege's account of judgment. 

\subsection{Frege's account of judgment} \label{subsec:judg}

Given that generality is expressed in a judgment, a few words on Frege's account of judgment may be in order. First of all, let me emphasize that all axioms and theorems of the ideography take the form of judgments that state that a thought~is~true. For Frege, a judgment is the acknowledgment of the truth of a thought, or, equivalently, the acknowledgment that the proposition expressed by a sentence refers to the True: 

\begin{quote}
	with `$2 + 3 = 5$' only a truth-value is designated, without its being said which
	one of the two it is. Moreover, if I wrote `$(2 + 3 = 5) = (2 = 2)$' and presupposed that
	one knows that $2 = 2$ is the True, even then I would not thereby have asserted that
	the sum of $2$ and $3$ is $5$; rather I would only have designated the truth-value of: that
	`$2 + 3 = 5$' refers to the same as `$2 = 2$'. We are therefore in need of another special
	sign in order to be able to assert something as true. To this end, I let the sign
	precede the name of the truth-value, in such a way that, e.g., in
	
	\vspace*{-1.3em}
	
	$$\text{`}\Fa 2^2 = 4\text{'}$$
	
	\vspace*{-1.3em}
	
	it is asserted that the square of $2$ is $4$. I distinguish the \textit{judgement} from the \textit{thought}
	in such a way that I understand by a \textit{judgement} the acknowledgement of the truth of
	a \textit{thought}. (\gga~\S5)
\end{quote}

I want to stress that ``acknowledgment'' is a success word. It is clear that for Frege a~judgment ``is not the mere grasping of a thought, but the admission of its truth''~\cite[fn.7]{frege1892sense}. 
Therefore, in Frege a judgment is always synonymous with a correct judgment. This reveals that the turnstile sign of \gga~has a double role of expressing a judgment while at the same time making the judgment expressed~\citep{hare1970meaning}. 
Formally speaking, a judgment in the ideography consists of a sentence  preceded of a turnstile sign:

\vspace*{-1.3em}

$$\Fa[1] a. $$

Since sentences are names of objects, note that Frege does not respect the now standard distinction between terms and formulas in modern predicate logic. \citet[\S2.2]{landini2012frege} maintains that this orthodox reading is incorrect because, under his interpretation, formulas can be formed by appending a turnstile sign to terms of the ideography. So, 
he sees `$\Fa a$' as a formula whenever `$a$' is a term of the formal language.\footnote{So does \citet[p.30]{klement2001frege}, except that Landini's view is purely syntactic. } 
I am afraid this interpretation neglects the madeness associated with Frege's judgment stroke. 
Landini thinks to convey correctness with Frege's turnstile sign is to conflate matters of syntax with those of semantics in the ideography. 
To be fair, Landini admits that Frege can never be found attaching the turnstile sign to what he believes is a false sentence, but counters that ``Frege only says that attaching a turnstile `aims' to make a judgment''~\cite[p.31]{landini2012frege}. 
Landini emphasizes Frege does express his false Basic Law V with a turnstile. So, in his view, a turnstile sign can flank not just names of the True but any term. To conclude, \citet[p.34]{landini2012frege} contrasts the use of the turnstile sign `$\vdash$' in modern logic with Frege's use of his turnstile sign `$\Fa$'. 
The former is attached to a formula of the language to form a thesis (axiom or theorem) only expressible in the metalanguage, while, under his view, the latter is attached to a term to form a formula in the object language. 
Implicit in Landini's defense of his formula interpretation of judgments is the premise that Frege's turnstile sign is part of the object language. 
But as `$\Fa a$' goes beyond the formal theory by saying that the thought expressed by `$a$' in the object language refers to the True, I see this as evidence that it is not part of the object language. If that is correct, there is no confusion between syntax and semantics, as Landini worries. The fact that Frege used a turnstile sign to express Basic Law V in \gga~simply means, as \citet[p.8]{whitehead1997principia} point out, that he can be convicted of error. 

Moreover, if the madeness associated with the judgment stroke is denied, and Frege just meant to say that attaching a turnstile sign to a name ``aims'' to make a judgment, as Landini maintains, one may wonder how similar the roles of the turnstile in \gga~and the horizontal in \bs~would be. 
Expressing that something can become a judgment was precisely the original role of the horizontal in \bs~\S2. In the formal system of \bs~\text{} $\F A$ expresses that ``the content $A$ is judgeable'', which, from a purely syntactic perspective, is to say it is a formula~\cite[\S2]{bentzen2020perspectiva}. 
Furthermore, given that all Frege has to say about the madeness of judgments in \gga~\S5 and~\bs~\S2 amounts to essentially the same, it seems to me that for a conclusive vindication of the formula interpretation of judgments one has to explain how it can tell the judgments and judgeable contents (formulas) of \bs~apart. 
One may respond that the formula interpretation of judgments only applies to \gga~because in \bs~the judgment stroke already flanks a formula while in \gga~it flanks a term. But that would be to overlook that Frege sees a sentence as a name of a truth value in \gga.

Back to the discussion of generality in \gga~\S17. Frege states that a Roman marker can be transformed into a universally quantified sentence, thereby explaining what seems to be the only purpose of Roman markers in the ideography, formally speaking. This transition is illustrated as a rule of universal generalization that allows us to pass from a judgment presented with a Roman marker of a truth value to a judgment presented with a sentence: \\


\begin{equation} \label{eq:gen17}
\frac{\qquad\qquad\Fa[1] \Phi(x)\qquad}{\Faquant {a} \Phi(\mathfrak{a}).}
\end{equation} 

\vspace*{-1.em}

Like any other rules of inference in the ideography, universal generalization operates on judgments as the premises and conclusions of the rule.\footnote{
	It is worth noting that the ideography lacks a rule of conditional proof, which explains why Frege's universal generalization rule can be formulated without any limitations~\citep[p.33]{landini2012frege}.} 
Therefore, Frege implies that we have a judgment in the premise of this rule. Despite Frege's implicit assumption, a second look reveals that the account of judgment outlined in \gga~\S5 is only applicable in the context of sentences, not Roman markers. How can we acknowledge the truth of a thought when a Roman marker, not being a sentence, neither express a thought nor refers? 
It is clear that a ``judgment'' of this kind has to be understood in a different way. 
What we can acknowledge here is that no matter how the Roman letters are replaced with names the resulting name always refers to the True. 
It is not an acknowledgment of a true thought but of true thoughts whenever appropriate names are given. I will argue in \Cref{openterms} that this is the explanation that Frege should have adopted for such judgments. It will be seen below, however, that a different account is hinted at in \gga~\S32. 

Let me introduce some terminology for the sake of convenience. 
If we think of Roman markers as open terms and names as closed terms following~\cite{landini2012frege}, but use simple type theory as the framework of our investigation instead of the predicate calculus, then we can see that judgments in the ideography are also either categorical or general hypothetical. So, when referring to the ideography in what follows I call a judgment ``categorical'' if it is composed of a turnstile sign followed by a name of a truth value or ``general hypothetical'' if a turnstile precedes a Roman marker that indicates a truth value. To give a concrete example, these are categorical and hypothetical judgments of the ideography, respectively: 

$$ \Facontent \spiritusasper{\varepsilon}(\varepsilon = \varepsilon ) \qquad\qquad ~ \Fa\Fconditional{\Fcontent f(x)=f(y)}{\Fcontent x=y.}$$

Frege is well aware of the subtle distinction involved here and seems to use it quite skillfully. It is explicitly stated, for instance, in \gga~\S26 and \S32, that a judgment in the ideography can only take the categorical and general hypothetical forms described above. Frege opens \gga~\S17 with a remark that Barbara cannot be derived from his previous rules of inference when the premises and conclusion are understood categorically in terms of universally quantified sentences, but only when represented general-hypothetically. 

Since the account of judgment given in \gga~\S5 is not applicable to general hypothetical judgments, Frege attempts to justify their presence in the ideography in an alternative way in \gga~\S32. Therein it is observed that there are in fact only categorical judgments in the ideography as long as we consider that with the rule of universal generalization~(\ref{eq:gen17}) a Roman marker can be transformed into a general sentence:

\begin{quote}
	Now, a concept-script
	proposition consists of a judgement-stroke with a name, or a Roman marker, of a
	truth-value. However, such a marker is transformed into the name of a truth-value
	when German letters are introduced for the Roman ones with concavities put in front,
	in accordance with §17. If we suppose this has been carried out, then there is only
	the case where the proposition is composed of the judgement-stroke and a name of a
	truth-value. (\gga~\S32)
\end{quote}

It is unclear how to connect this passage with the view exposed in \gga~\S5 in terms of acknowledgments, but we may presume that with a general hypothetical judgment $~\Fa[1] \Phi(x)$ we acknowledge, although only indirectly, the truth of $\Fquant {a} \Phi(\mathfrak{a})$ instead. One problem with this account will be discussed next in \Cref{openterms}: Frege appears to incur circularity. 

\section{Functionality and generality in the ideography} \label{openterms}

What is the rationale behind Frege's stipulation that only function names but not Roman markers are allowed to refer as incomplete expressions? What criterion does Frege have in the expository language of \gga~to tell them apart? They certainly differ in that the argument places of a function name are marked with small Greek letters while the incompleteness of a Roman marker is expressed with Roman letters. 
But to adopt this minor syntactic convention as the criterion for determining that function name may refer but Roman markers may only indicate would be unacceptable. 
It could be argued that function names only occur in the metalanguage but Roman markers are part of the object language. 
However, the question again arises as to what criteria are used to support his decision of doing so.
If Frege's distinction between them is meaningful, a satisfactory explanation cannot simply rely on his arbitrary convention of using small Greek letters as metavariables. 

To my mind, the problem is that both small Greek and Roman letters are attributed similar semantic roles as signs that ``only indicate indeterminately''. How do we distinguish those two kinds of incomplete expressions, function names and Roman markers, without begging the question by appealing to the very stipulation in question that Roman markers do not refer like function names but merely indicate? Frege does not provide a direct answer for this question, but the following criteria for the distinction may be drawn from his expository remarks on function names and Roman markers already discussed in the two previous sections: 

\begin{itemize}
	\item \textit{functionality}: a function name expresses functionality as something that establishes a connection between the objects whose names we put for a small Greek letter and the objects that result as referents of the resulting name, as seen in \Cref{objectandfunction};
	
	\item \textit{generality}: a Roman marker expresses generality in a judgment with its Roman letters as something that can be transformed into a general sentence, as seen in \Cref{romanmarkers}. 
\end{itemize}

In this section I will argue that Frege lacks a clear criterion for distinguishing between function names and Roman markers because his views on functionality apply to Roman markers and the same can be said of his views on generality applying just as well to function names. 

First, let us note that the criterion of functionality laid down in \gga~\S1, and quoted and discussed in \Cref{objectandfunction}, applies to Roman markers. 
We have already seen in \Cref{romanletters} that a Roman marker can relate the objects whose names we put for a Roman letter and the objects that appear as referents of the substituted expression. If Roman markers are not functional just as function names are, then functionality cannot be manifested in the part of an incomplete expression that remains without the variables, here either marked with small Greek or Roman letters, both of which, according to Frege, only indicate indeterminately. 
Why is `$\xi = \zeta$' functional but not `$a = b$', for instance, considering that both expressions result in referential names when their variables are replaced with referential names? The parts of the expressions that remain without the variables are syntactically identical, so we can say with Frege that the ``nature of the function'' must ``lie'' in both incomplete expressions. That is, it must be admitted, if we wish to conform to Frege's criterion of functionality, that both function names and Roman markers are functional in \gga. 

Second, what can be said about the applicability of the generality criterion to function names rather than only Roman markers? Recall that the generality expressed by Roman markers is formally made precise with the transition (\ref{eq:gen17}) from \gga~\S17. 
Now consider a transition where now we have a passage from a function name to a general sentence, stating that given a function $\Phi(\xi)$ assumed to always have the True as a value for every argument, we can apply the second-level function $\Fquant {a} \phi(\mathfrak{a})$ to obtain a true thought $\Fquant {a} \Phi(\mathfrak{a})$:  

\begin{equation} \label{eq:conceptabstracttion}
\frac{\qquad\quad \Fa[1] \Phi(\xi)\qquad}{\Faquant {a} \Phi(\mathfrak{a}).}
\end{equation} 

Of course, this controversial transition is nowhere to be found in \gga~and it is completely ungrammatical from the standpoint of its object language. 
Readers will no doubt object that we cannot see this as a transition on the same level as (\ref{eq:gen17}) because in (\ref{eq:conceptabstracttion}) we do not actually have a judgment in the premise. 
But I believe it would be a mistake to think that we are in a better position to justify (\ref{eq:gen17}) as a transition. In its premise we have a Roman marker and $\Fa[1] \Phi(x)$ is not a categorical judgment either. For this reason, Frege cannot justify (\ref{eq:gen17}) unless he offers in advance an explanation of what makes its premise a judgment. 

This amounts to giving an account of general hypothetical judgments. 
But earlier in \Cref{romanmarkers} we saw Frege defending in \gga~\S32 the presence of general hypothetical judgments in the ideography by reducing them to categorical judgments using the very transition (\ref{eq:gen17}) that transforms Roman markers into general sentences. 
This is a circular strategy. It attempts to explain the introduction of general hypothetical judgments by taking for granted that (\ref{eq:gen17}) is a valid transition; yet this implies that both its conclusion and its premise, which happens to take a general hypothetical form, are judgments of the ideography. 

This is why I suggested in \Cref{romanmarkers} that it is more promising to understand a general hypothetical judgment as an acknowledgment that an incomplete expression expresses a truth thought whenever it is filled with referential names of suitable levels. 
Surely $~\Fa[1] \Phi(\xi)$ can also be understood in those terms if the function has the True as value for every argument. 
But even if one wants to stick to Frege's circular account of general hypothetical judgments I see no reason to accept (\ref{eq:gen17}) as a valid transition but not (\ref{eq:conceptabstracttion}), unless one begs the question by contending that Roman markers are not function names. 
Given that both Roman markers and function names seem equally transformable into general sentences, generality cannot be a reliable criterion to distinguish between Roman markers and function names.

In sum, the claim that function names refer but Roman markers fail to do so is unfounded, given that Frege does not present any convincing criterion to distinguish between Roman markers and function names. 
By insisting on distinguishing function names and Roman markers, Frege creates a puzzling distinction between two classes of incomplete expressions.  
This observation paves the way for the open term interpretation: no conflict should arise in the reading of both function names and Roman markers as open terms. 
It also reveals that Frege could not be justified in keeping function names outside the object language and Roman markers in the object language of \gga. 

\section{The open term interpretation}\label{simpletypetheory} 

The open term interpretation may be briefly described as the type-theoretic rendering of Frege's view of functions as incomplete entities. 
It embraces \citepos{church1939schroeder} reading of Frege's functions as incompleted abstractions and addresses~\citepos{klev2014categories} concern about cases of functions being the arguments to other
functions by giving an interpretation of second-level functions without function types. 
Now, the correlation between  function names in 
\gga~and open terms in simple type theory should be clear at this point. 
Following the functionality criterion laid down in \Cref{openterms}, a function name in \gga~must establish a connection between the objects whose names are put for a variable and the objects that are values of the complete expression. But, following \cite{martinlof1982constructive}, this is exactly what it means for an open term with explicit occurrences of variables to have a type. 

Thus, if we were to view the primitive first-level function names of \gga~as terms of a simple type theory, they would have to be open terms of a ground type and not closed terms of the function type. The resulting formalism would be a fragment of simple type theory that lacks a function type and the weakening property, but includes the following typing judgments, where $\iota$ is a type of individuals and $o$ a type of truth values:

$$\vdash \hypjdg{x : \iota}{\F[1] x : o}, \quad\quad \vdash \hypjdg{x : \iota}{\Fncontent x : o}, \quad\quad \vdash \hypjdg{x, y : \iota}{\Fcontent\Fconditional{\Fcontent y}{\Fcontent x} : o}, $$ 

$$\vdash \hypjdg{x, y : \iota}{x = y : o}, \quad\quad \vdash \hypjdg{x : \iota}{\backslash x : \iota}.$$

\noindent Moreover, because every truth value is an object, we require that:

$$\frac{\vdash \hypjdg{\Gamma}{a : o}}{\vdash \hypjdg{\Gamma}{a : \iota}.}$$

Here open terms inhabit ground types and first-level functions are open terms whose argument places are designated by free typed variables. We must reject weakening because it goes against Frege's conviction that a function name must have explicit argument places. If we had this property it would not be possible to view open terms as function names. Given a closed term such as $\spirituslenis{x}(\F[1]x) : \iota$ (a term for a value-range now written in simplified notation without small Greek letters), we could obtain, by weakening, an open term $\hypjdg{y : \iota}{\spirituslenis{x}(\F[1]x) : \iota}$ that would otherwise be a function name of one argument, except that no argument places actually occur in the expression, in contrary to Frege's view. 

Function application takes the form of substitution (\ref{eq:ttsubst}) of open terms. One benefit of this view is that application for functions of two arguments is handled in a way that resembles the method of double completion introduced in \gga~\S4. From $\hypjdg{x, y : \iota}{f(x,y) : \iota}$ and $a : \iota$ we may obtain either $\hypjdg{x : \iota}{f(x,a) : \iota}$ or $\hypjdg{y : \iota}{f(a,y) : \iota}$. Only after yet another completion, say, with $b : \iota$, we arrive at a closed term  $f(a,b) : \iota$. Being a closed term, this is called the value of the function. 
Just as with Frege, there is no fixed application order for functions of two arguments, but we only arrive at a value after a double completion. 
As noted before, Curried functions fail to represent Frege's view in both respects. 

Finally, it should be mentioned that  under this open term interpretation the definition of new functions (in the old-fashioned sense) can be done by uniformly substituting complex function signs for object language function variables. This corresponds to \citepos{landini2012frege} interpretation of comprehension of functions in \gga. 
Substitution is made possible with the introduction of the following rule of inference, where $a$ is free for $x$ in $f$:

$$ \frac{\vdash \hypjdg{\Gamma}{f [a/x] : \tau} \quad \vdash a : \sigma}{\vdash \hypjdg{\Gamma, x : \sigma}{f : \tau}.}$$

\subsection{Second-level function names} \label{secondlevel}

The view of first-level functions as open terms bridges the gap between the approaches to function abstraction in the ideography and simple type theory. In both systems, the operation may be viewed as taking an open term depending on a distinguished free variable and resulting in a closed term.  The main difference is that for Frege the abstracted term is in the ground type of individuals and not in a function type, unlike in (\ref{eq:lambdaabst}). 
So instead we have: 

\begin{equation} \label{eq:abstraction}
\frac{\vdash \hypjdg{x : \iota}{a : \iota}}{ \vdash \spirituslenis{x}a : \iota.}
\end{equation}

\newpage

However, the reading of Frege's function names as open terms becomes less simple whenever higher order functions are concerned due to the absence of function types. Open terms are formed by means of general hypothetical judgments, and those cannot be part of other general hypothetical judgments in the same way that higher-level functions may take lower-level functions as arguments. 
Fortunately, there is a very natural way to circumvent our limitation and extend our interpretation to second-level functions. I have just explained the value-range operator as an inference rule that takes an open term in the type of individuals $\iota$ to a closed term in that type. The same can be said of first-order quantification using the type of truth values $o$. I now use modern notation for the sake of convenience:

\begin{equation} \label{eq:quantifier}
\frac{\vdash \hypjdg{x : \iota}{f : o}}{ \vdash \forall (x : \iota) f : o.} 
\end{equation}

Another way in which this is a natural interpretation is the following. Following the observation from \Cref{romanmarkers} that Roman markers behave as open terms, (\ref{eq:quantifier}) also interprets the transition  (\ref{eq:gen17}) from Roman markers to general sentences proposed in \gga~\S17, revealing that Frege's distinction between Roman markers and function names is ultimately superfluous when seen through the lenses of simple type theory. 
Both are open terms that can be transformed into closed terms representing value-range or quantified propositions. 

This interpretation sheds light on a noteworthy relation between function abstraction and quantification. 
Under a correspondence famously discovered by \cite{howard1980formulae} between propositions and types, the terms of the function type can be regarded as proofs of implications. Just as simple type theory can be viewed as a generalized form of propositional logic, a more sophisticated framework known as ``constructive type theory''~\citep{martinlof1975intuitionistic} generalizes predicate calculus. 
The theory features dependent types, a constructive extension of Frege's conception of predicates as functions that assign objects to reified truth values. 
From a suitable open term of a type we can always either obtain a value-range term from (\ref{eq:abstraction}) or a quantified term from (\ref{eq:quantifier}), but never both in a single rule. In dependent type theory, this conflict is resolved by keeping function abstraction and quantification apart with one rule that at the same time emphasizes both the functionality and generality of an open term: 

\begin{equation} \label{eq:dependentfunc}
\frac{\vdash \hypjdg{\Gamma, x : \sigma}{a : f}}{\vdash \hypjdg{\Gamma}{\fun x a : \forall(x : \sigma) f} .}
\end{equation}

Under the propositions-as-types correspondence, quantification $\forall(x : \sigma) f$ can be seen as a generalized type of functions $\funtype{(x : \sigma)}{f}$ whose values for any $a : \sigma$ have type $f[a/x]$. 
Thus, $f$ is an indexed family of types associating every term of a base type to another type. 
This also means that instead of $\hypjdg{x : \iota}{f(x) : o}$, in constructive type theory we require $\hypjdg{x : \iota}{f(x) : \mathcal{U}}$, where $\mathcal{U}$ is a universe type, a type whose terms are themselves types. 
Of~course, I am not suggesting that Frege would be prepared to accept~(\ref{eq:dependentfunc}). I doubt he could even make sense of it as a rule of inference in his wildest dreams. All I wish to point out is that, since lambda terms are function terms, Frege's struggle with the distinction between function names and Roman markers, which are both here open terms, is finally settled in constructive type theory by keeping functionality on the level of terms and generality on the level of types. 

\subsection{Third-level function names} \label{thidlevel}

One objection to the open term interpretation is that is unclear how the approach to second-level functions could be generalized to third-level functions because we obviously do not have something that operates on rules of inferences in type theory. Such an interpretation would therefore have to be studied at a metainferential level, a task which I am not willing to undertake here since I am not confident about its philosophical significance. 
The purpose of this paper is not to provide an entirely accurate representation of Frege's theory of function levels assuming simple type theory in the background. Instead, the objective is to analyze and compare the function type and open term interpretations in order to determine which one aligns more closely with Frege's theory of function levels.

\section{Concluding remarks} \label{conclusion}

I have discussed two representations of Frege's theory of function levels in simple type theory, namely, the function type and open term interpretations. The former states that a type of functions is presupposed at the level of the expository language of \gga. While this interpretation is very influential in the literature, it leaves out function names of two arguments and goes against Frege's view of function names as incomplete expressions. 
On the other hand, Frege's lack of criterion to distinguish function names from Roman markers and his account of functionality seem to favor the open term interpretation. Thus understood, \gga~only presupposes a ground type of individuals $\iota$ and truth-values $o$ assigned to every closed term in the language by virtue of the fact that they are complete expressions. 
There are technical complications in extending the open term interpretation to functions of third level, but it brings Frege's and Church's theories of abstraction together under one common framework for open terms and sheds new light on the type-theoretic rendering of the view that function names are incomplete expressions. 


\bibliographystyle{plainnat}
{\linespread{0.7}\selectfont\bibliography{ref}}

\vspace*{-1em}

\begin{acknowledgments}
	I would like to thank Ansten Klev, Marco Panza, Marco Ruffino, Will Stafford, Göran Sundholm, Wim Vanrie, and two anonymous reviewers for their helpful comments on early drafts of this paper. Most of the material here has been presented in different online seminars and conferences over the past three years, in particular the {OCIE Seminar in the History and Philosophy of Mathematics and Logic}, {Logicians in Quarantine}, {I Colóquio NULFA-NuLFiC de Lógica e Filosofia da Matemática}, and {Ch\'a das 5}. I also thank Alessandro Duarte, Erich Reck, and Edward Zalta for stimulating discussions. This research was partly supported by the AFOSR grant FA9550-18-1-0120 and the Lumina quaeruntur fellowship number LQ300092101 from the Czech Academy of Sciences. 
\end{acknowledgments}

\vspace*{-2em}

\address

\end{document}